\documentclass{article}
\usepackage{amsmath,amsthm,amssymb,amscd,fancybox,epsfig,float,subfigure}
\usepackage{graphics}
\pdfoutput=1

\newtheorem{definition}{Definition}[section]
\newtheorem{remark}{Remark}[section]

\newcommand{\beqn}{\begin{equation}}
\newcommand{\eeqn}{\end{equation}}

\def\xb{{\mathbf x}}
\def\yb{{\mathbf y}}

\def\E{{\mathbf E}}
\def\H{{\mathbf H}}
\def\bE{{\mathbf E}}
\def\bH{{\mathbf H}}

\def\A{{\mathbf A}}

\def\J{{\mathbf J}}

\def\K{{\mathbf K}}

\def\n{{\mathbf n}}

\newcommand{\by}{\ensuremath{{\bf{y}}}}

\newcommand{\rthree}{\ensuremath{{\mathbb{R}}^3}}

\newcommand{\ben}{\begin{equation*}}
\newcommand{\een}{\end{equation*}}
\newcommand{\be} [1] {\begin{equation} \label{#1}}
\newcommand{\ee}{\end{equation}}

\newtheorem{lemma}{Lemma}

\title{Fast multi-particle scattering: a hybrid solver for the Maxwell equations in 
microstructured materials}

\author{Z. Gimbutas
and L. Greengard
\thanks{Courant Institute of Mathematical Sciences, 
         New York University, 251 Mercer Street,
         New York, NY 10012-1110.
{{\em email}: {\sf {gimbutas@courant.nyu.edu}},
{\sf {greengar@courant.nyu.edu}}}. This work was supported by
the U.S. Department of Energy 
under contract DEFGO288ER25053 and by 
the Air Force Office of Scientific Research under MURI grant
FA9550-06-1-0337 and under NSSEFF Program Award FA9550-10-1-0180 
.}
}

\begin{document}
\maketitle

\begin{abstract}
A variety of problems in device and materials design require the 
rapid forward modeling of Maxwell's equations in complex micro-structured
materials. By combining high-order accurate integral equation methods with
classical multiple scattering theory, we have created an effective simulation
tool for  materials consisting of an isotropic background in which are
dispersed a large number of micro- or nano-scale metallic or dielectric inclusions.
\\
\noindent {\bf Keywords}: Maxwell equations, multiple scattering, meta-materials,
fast multipole method

\end{abstract}

\section{Introduction}

\begin{figure}
\centering
\includegraphics[width=5in]{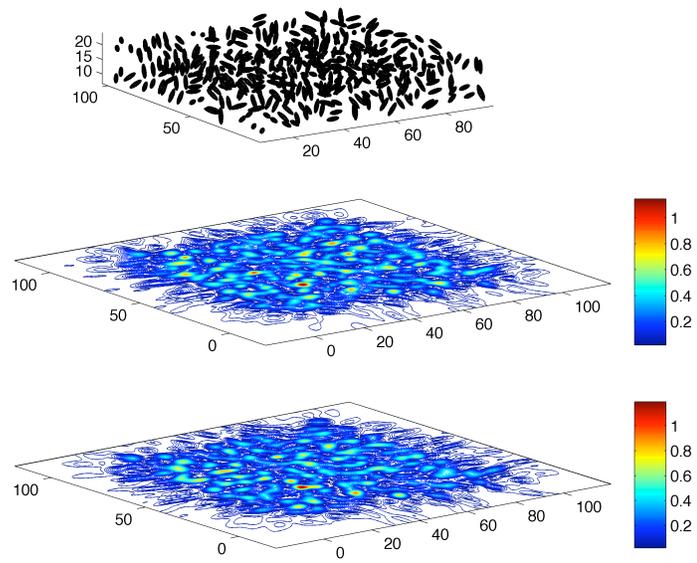} 
\caption{\sf Two hundred gold ellipsoid pairs
are randomly oriented in the region $[0,100] \times [0,100] \times [0,20]$ and illuminated from above by a plane wave in TE polarization. 
The transmitted $z$-component of the Poynting vector is plotted
on planes at $z=-4$ and $z=-8$. The wavelength is $2 \pi$ so that
the particles are approximately one wavelength
in size, and the region is about $15 \times 15 \times 3$ wavelengths is size. 
\label{f:1}}
\end{figure}

We describe in this paper a simulation method for Maxwell's equations suitable for 
microstructured materials consisting of separated inclusions which are embedded in a homogeneous
background (Fig. \ref{f:1}).
In practice, it is often the case that the shape and permittivity of the inclusions
are fixed and that one seeks to optimize their placement to create a specific electromagnetic
response. Each new configuration, however, requires the solution of the full Maxwell equations.
If there are thousands of inclusions in an electrically large region (many wavelengths in size), the
calculation is generally too expensive to carry out within a design loop. 

In oder to accelerate such calculations, we have coupled complex geometry Maxwell solvers with
multiple scattering theory. Using the hybrid solver, 
calculations such as the one depicted in Fig. \ref{f:1} require
only a few minutes on a single CPU, despite the fact that there are a million degrees of freedom
needed to describe the full geometry (and there would be orders of magnitude more points
needed in a finite difference or finite element discretization).

Our method, which we refer to as fast multi-particle scattering (FMPS), is based on a two step procedure. 
First, we enclose a representative scatterer, such as a single pair of gold nanorods, in 
a sphere $S$. We then build the scattering matrix for this nano-structure (described below) using
integral equation techniques. 
The solution to the full Maxwell equations can then be obtained in geometries with $N$
inclusions ($N = 200$ in Fig. \ref{f:1}), by solving the multiple-scattering problem where
the inclusions have been replaced by their scattering matrices. Not only does this reduce
the number of degrees of freedom required, but we have effectively precomputed
the solution operator for each inclusion in isolation, so that the linear system we solve
by iteration on the multi-sphere system
is well-conditioned. Further, the fast multipole method (FMM) reduces the cost of each iteration 
from $O(N^2)$ to $O(N \log N)$ and is particularly efficient when applied to this problem.

The principal limitations of the method are (1) that  some modest separation
distance between inclusions is required and (2) that some of 
the efficiency is based on the fact that only a few distinct nanoparticle types are allowed. 
In many experimental settings, both conditions are satisfied. We will return to a discussion
of these limitations in our concluding remarks.

\section{Maxwell's equations and the Debye-Lorenz-Mie formalism}

Working in the frequency domain and assuming a time dependence 
of $e^{-i\omega t}$, Maxwell's equations in a linear, isotropic material take the form
\begin{align}
  \nabla \times {\bf H}^{tot} &= -i\omega \epsilon \,  {\bf E}^{tot},
\label{meq1} \\
  \nabla \times {\bf E}^{tot} &= i\omega \mu {\bf H}^{tot},  \nonumber
\end{align}
where ${\bf E}^{tot}$ and ${\bf H}^{tot}$ are the total electric and magnetic
fields.  $\epsilon$ is the permittivity of the medium and $\mu$ its permeability.
We are mainly interested in dielectric inclusions embedded in a background medium,
but will consider perfect conductors briefly  at the end of this section.
The total fields $({\bf E}^{tot},{\bf H}^{tot})$ can be written as the
sums of the incident fields $({\bf E}^{in},{\bf H}^{in})$, defined
only in the exterior region, and scattered fields $({\bf E},{\bf H})$ defined in both the inclusions and the
exterior:
\begin{eqnarray}
  {\bf E}^{tot} &=& {\bf E}^{in} + {\bf E}, \nonumber \\
  {\bf H}^{tot} &=& {\bf H}^{in} + {\bf H}  \label{fbgvi3a} .
\end{eqnarray}
It is well-known \cite{Jackson, Papas} that at dielectric
interfaces, the Maxwell equations (\ref{meq1}) are uniquely solvable when supplemented by the 
the continuity conditions: 
\begin{eqnarray}
 \left[  {\bf n} \times {\bf E}^{tot}  \right] &=& {\bf 0} \quad \Rightarrow \quad
 \left[ {\bf n} \times {\bf E} \right] = -
 \left[ {\bf n} \times {\bf E}^{in} \right]
    \nonumber \\
   \left[  {\bf n} \times {\bf H}^{tot} \right] &=& {\bf 0} \quad \Rightarrow \quad
 \left[ {\bf n} \times {\bf H} \right] = -
 \left[ {\bf n} \times {\bf H}^{in} \right]
\label{mwaveint1} 
\end{eqnarray}
and the Silver-M{\" u}ller radiation conditions on the scattered field. The expression 
$[ {\bf n} \times {\bf F} ]$  is used to denote the jump in the tangential components of the 
vector field ${\bf F}$ at a point on the interface.

\subsection{Debye Potentials}

About a century ago, Debye, Lorenz, and Mie \cite{Debye,Lorenz,Mie} independently solved
the problem of scattering from a single sphere by using separation of variables.
Without entering into the derivation, it is straightforward to verify that
\begin{eqnarray} 
\bE(\xb) &=& \nabla \times \nabla \times (\xb v(\xb) 
\, + \, i \omega \epsilon \nabla \times (\xb u(\xb)) \nonumber \\
\bH(\xb) &=& \nabla \times \nabla \times (\xb u(\xb)) 
- i \omega \mu \nabla \times (\xb v(\xb)) 
\label{debyerep}
\end{eqnarray}
represent an electromagnetic field, where $\xb$ denotes the position vector with respect to the 
sphere center, so long as the 
{\em Debye potentials} $u,v$ satisfy the scalar Helmholtz equation 
\[ \Delta u + k^2 u = 0, \ \Delta v + k^2 v = 0\ , \]
with Helmholtz parameter (wave number) $k^2 = \omega^2 \epsilon \mu$. 
In the exterior of a sphere, the Debye potentials $u,v$ can be represented by the multipole expansions
\begin{eqnarray}
 u^{ext}(r,\theta,\phi) &=& \sum_{n=0}^\infty \sum_{m=-n}^n b_{n,m} h_n(kr) Y_n^m(\theta,\phi) \nonumber \\
 v^{ext}(r,\theta,\phi) &=& \sum_{n=0}^\infty \sum_{m=-n}^n a_{n,m} h_n(kr) Y_n^m(\theta,\phi)  
\label{eq:uvext}
\end{eqnarray}
where $(r,\theta,\phi)$ are the spherical coordinates of the point $\xb$ with respect to the sphere center, 
 $h_n(r)$ is the spherical Hankel function of order $n$, and
$Y_n^m(\theta,\phi)$ is the usual spherical harmonic of order $n$ and degree $m$. 
The resulting electromagnetic field then also satisfies the appropriate radiation conditions at 
infinity.
In the interior of a sphere,  $u$ and $v$ can be represented by the local expansions
\begin{eqnarray}
u^{int}(r,\theta,\phi) &=& \sum_{n=0}^\infty \sum_{m=-n}^n d_{n,m} j_n(kr) Y_n^m(\theta,\phi) \nonumber \\
v^{int}(r,\theta,\phi) &=& \sum_{n=0}^\infty \sum_{m=-n}^n c_{n,m} j_n(kr) Y_n^m(\theta,\phi)  
\label{eq:uvint}
\end{eqnarray}
where $j_n(x)$ is the spherical Bessel function of order $n$. 

\begin{remark}
To improve readability, we will abbreviate 
\[    \sum_{n=0}^\infty \sum_{m=-n}^n     \qquad {\rm as} \qquad  \sum_{n,m}  \]
and the truncated sum
\[    \sum_{n=0}^p \sum_{m=-n}^n     \qquad {\rm as}  \qquad \sum^p_{n,m} \, . \]
It is straightforward to verify that the total number of terms in the truncated summation
is $(p+1)^2$.
\end{remark}

\subsection{Single sphere scattering} \label{sec:single}

Suppose now that one is interested in scattering from a single dielectric sphere $S$
of radius $R$ with permittivity $\epsilon_1$, permeability $\mu_1$, and Helmholtz parameter $k_1 = 
\sqrt{\omega^2 \epsilon_1 \mu_1}$, in response to an incoming 
field $(\bE^{in},\bH^{in})$. The external medium is assume to have 
permittivity $\epsilon_0$, permeability $\mu_0$, and Helmholtz parameter $k_0 = 
\sqrt{\omega^2 \epsilon_0 \mu_0}$. Then the scattered field can be represented 
by (\ref{debyerep}) with $k = k_1$ in (\ref{eq:uvint})
for $(r,\theta,\phi)$ inside $S$ and 
by (\ref{debyerep}) with $k = k_0$ in (\ref{eq:uvext})
for $(r,\theta,\phi)$ outside $S$. 

Let us denote by $\bE_0, \bH_0$ the scattered field in the 
exterior domain and by  $\bE_1,\bH_1$ the scattered field inside $S$. Then
\begin{eqnarray*} 
\bE_0(\xb) &=&  \sum_{n,m} a_{n,m} 
\nabla \times \nabla \times (\xb \phi_{n,m}^{k_0}) 
\, + \, i \omega \mu_0 \sum_{n,m} b_{n,m} \nabla \times (\xb \phi_{n,m}^{k_0}) \\
\bH_0(\xb) &=&  \sum_{n,m} b_{n,m} 
\nabla \times \nabla \times (\xb \phi_{n,m}^{k_0}) 
\, - \, i \omega \epsilon_0 \sum_{n,m} a_{n,m} \nabla \times (\xb \phi_{n,m}^{k_0}) 
\end{eqnarray*}
where 
$ \phi_{n,m}^k(\xb) =  \phi_{n,m}^k \left[ r,\theta,\phi \right] = h_n(k r) Y_n^m(\theta,\phi) $ and
\begin{eqnarray*} 
\bE_1(\xb) &=&  \sum_{n,m} c_{n,m} 
\nabla \times \nabla \times (\xb \psi_{n,m}^{k_1}) 
\, + \, i \omega \mu_1 \sum_{n,m} d_{n,m} \nabla \times (\xb \psi_{n,m}^{k_1}) \\
\bH_1(\xb) &=&  \sum_{n,m} d_{n,m} 
\nabla \times \nabla \times (\xb \psi_{n,m}^{k_1}) 
\, - \, i \omega \epsilon_1 \sum_{n,m} c_{n,m} \nabla \times (\xb \psi_{n,m}^{k_1}) 
\end{eqnarray*}
where 
$\psi_{n,m}^{k}(\xb) =  \psi_{n,m}^{k} \left[ r,\theta,\phi \right] = j_n(k r) Y_n^m(\theta,\phi). $

We may also expand
$(\bE^{in},\bH^{in})$ in terms of spherical harmonics on the surface of $S$:
\begin{eqnarray*} 
\bE^{in}(\xb) &=&  \sum_{n,m} \alpha_{n,m}
\nabla \times \nabla \times (\xb \psi_{n,m}^{k_0}) 
\, + \, i \omega \mu_0 \sum_{n,m} \beta_{n,m} \nabla \times (\xb \psi_{n,m}^{k_0}) \\
\bH^{in}(\xb) &=&  \sum_{n,m} \beta_{n,m}
\nabla \times \nabla \times (\xb \psi_{n,m}^{k_0}) 
\, - \, i \omega \epsilon_0 \sum_{n,m} \alpha_{n,m} \nabla \times (\xb \psi_{n,m}^{k_0}) 
\end{eqnarray*}
 
All of the spherical
harmonic modes uncouple for fixed $n,m$, allowing for the determination of
$(a_{n,m}, b_{n,m}, c_{n,m}, d_{n,m})$ from the data
$(\alpha_{n,m}, \beta_{n,m})$
by applying the interface conditions (\ref{mwaveint1}). 
After some algebra (see, for example, \cite{Born,MTL}), one obtains two uncoupled
linear systems of the form
\begin{eqnarray} 
\left( \begin{array}{cc}
H_n(k_0R) & - J_n(k_1R) \\
\epsilon_0 h_n(k_0R)  & - \epsilon_1 j_n(k_1R) 
\end{array} \right) 
\left( \begin{array}{c} a_{n,m} \\ c_{n,m} \end{array} \right)  
= 
 \left( \begin{array}{c}  -J_n(k_0R) \alpha_{n,m} \\
 -\epsilon_0 j_n(k_0R) \alpha_{n,m}
 \end{array} \right) 
\label{Mie1}  \\
\left( \begin{array}{cc}
H_n(k_0R) & - J_n(k_1R) \\
\mu_0 h_n(k_0R)  & - \mu_1 j_n(k_1R) 
\end{array} \right) 
\left( \begin{array}{c}
 b_{n,m} \\
 d_{n,m}
 \end{array} \right)   = 
 \left( \begin{array}{c}
  -J_n(k_0R) \beta_{n,m} \\
 -\mu_0 j_n(k_0R) \beta_{n,m}
 \end{array} \right) 
\label{Mie2}
\end{eqnarray}
where $H_n(z) = [ h_n(z) + z h_n'(z)]$,  $J_n(z) = [ j_n(z) + z j_n'(z)]$. 

\begin{definition} \label{scatmatdef}
The mapping from incoming coefficients  $(\alpha_{n,m}, \beta_{n,m})$ to the 
outgoing coefficients $(a_{n,m},b_{n,m})$ is referred to as the {\em scattering
matrix} and denoted by {\cal S}.
\end{definition}

\subsection{Perfect conductors}

If the sphere $S$ is a perfect conductor, the corresponding boundary conditions are that the 
tangential components of the total electric field are zero \cite{Jackson, Papas}:
\begin{equation}
{\bf n} \times {\bf E}^{tot}  = {\bf 0} \quad \Rightarrow \quad
{\bf n} \times {\bf E}  = -  {\bf n} \times {\bf E}^{in} \, . 
\label{pec1} 
\end{equation}
In that case, the interior field is identically zero and the 
scattered matrix is given by
\begin{align} 
 a_{n,m}  &=  -(J_n(k_0R)/H_n(k_0R)) \alpha_{n,m} \nonumber \\
 b_{n,m} &=   -(j_n(k_0R)/h_n(k_0R)) \beta_{n,m}
\label{pecsolve}  
\end{align}

\section{Scattering from multiple spheres} \label{multspherescat}

Suppose now that one is interested in scattering from $M$ disjoint dielectric spheres,
where each sphere $S_l$ has  radius $R_l$ and  $k_l = 
\sqrt{\omega^2 \epsilon_l \mu_l}$. The external medium and incoming field are as above.
Then, the incoming field can be represented on the surface of $S_l$ by the expansion
\begin{eqnarray*} 
\bE_l^{in} =  \hspace{-.2in} && \sum_{n,m} \alpha_{n,m}^{l} 
\nabla \times \nabla \times (\xb_l \psi_{n,m}^{k_0}(\xb_l) ) 
\, + \, i \omega \mu_0 \sum_{n,m} \beta_{n,m}^{l} \nabla \times (\xb_l \psi_{n,m}^{k_0}(\xb_l) ) \\
\bH_l^{in} = \hspace{-.2in} &&  \sum_{n,m} \beta_{n,m}^{l} 
\nabla \times \nabla \times (\xb_l \psi_{n,m}^{k_0}(\xb_l) ) 
\, - \, i \omega \epsilon_0 \sum_{n,m} \alpha_{n,m}^{l} \nabla \times (\xb_l \psi_{n,m}^{k_0}(\xb_l) ) ,
\end{eqnarray*}
while the scattered field in the interior of $S_l$ can be represented by the 
expansion
\begin{eqnarray} 
\bE_l = \hspace{-.2in} && \sum_{n,m} c_{n,m}^{l} 
\nabla \times \nabla \times (\xb_l \psi_{n,m}^{k_l}(\xb_l) ) 
\, + \, i \omega \mu_l \sum_{n,m} d_{n,m}^{l} \nabla \times (\xb_l \psi_{n,m}^{k_l}(\xb_l) ) 
\hspace{.3in}  \label{Elint} \\
\bH_l = \hspace{-.2in} &&  \sum_{n,m} d_{n,m}^{l} 
\nabla \times \nabla \times (\xb_l \psi_{n,m}^{k_l}(\xb_l) ) 
\, - \, i \omega \epsilon_l \sum_{n,m} c_{n,m}^{l} \nabla \times (\xb_l \psi_{n,m}^{k_l}(\xb_l) ) .
\hspace{.3in} \label{Hlint}
\end{eqnarray}
Here, $\psi_{n,m}^{k}(\xb_l)= j_n(k r_l) Y_n^m(\theta_l,\phi_l)$
is computed in terms of the spherical coordinates $(r_l,\theta_l,\phi_l)$ of a point $\xb_l$
with respect to the center of  $S_l$.

The scattered field in the exterior of all the spheres can be represented by a sum of
outgoing expansions, one centered on each sphere.
\begin{eqnarray*} 
\bE_0 = \sum_{l=1}^M \sum_{n,m} a_{n,m}^{l} 
\nabla \times \nabla \times (\xb_l \phi_{n,m}^{k_0}(\xb_l) ) 
\, + \, i \omega \mu_0  \sum_{l=1}^M \sum_{n,m} b_{n,m}^{l} \nabla \times (\xb_l \phi_{n,m}^{k_0}(\xb_l)  ) \\
\bH_0 =  \sum_{l=1}^M  \sum_{n,m} b_{n,m}^{l} 
\nabla \times \nabla \times (\xb_l \phi_{n,m}^{k_0}(\xb_l)  ) 
\, - \, i \omega \epsilon_0 \sum_{l=1}^M  \sum_{n,m} a_{n,m}^{l} \nabla \times (\xb_l \phi_{n,m}^{k_0}(\xb_l)  ) .
\end{eqnarray*}
For a point $\xb$ exterior to all spheres, the function  
$\phi_{n,m}^{k_0}(\xb_l) \equiv  h_n(k_0 r_l) Y_n^m(\theta_l,\phi_l)$, where $\xb_l = (r_l,\theta_l,\phi_l)$, the latter being the spherical coordinates of $\xb$ with respect to the center of  $S_l$. The coefficients 
$(a_{n,m}^l, b_{n,m}^l, c_{n,m}^l, d_{n,m}^l)$ are all unknowns. They are determined by a linear system
that imposes the dielectric interface condition (\ref{mwaveint1}) on each sphere boundary.
Unlike the case of a single sphere, however, it is no longer trivial to solve for these unknowns, since the incoming field experienced on each sphere is due, not only to the known incoming field 
$(\bE^{in},\bH^{in})$,
but to the field scattered by all the other spheres. This results in a dense linear system
involving all of the unknowns, whose solution accounts for all of these {\em multiple scattering}
interactions. 

\subsection{Translation operators for multiple scattering} \label{sec:tml}

Fortunately, the outgoing Debye expansion on sphere $S_j$ can be analytically converted to an
incoming expansion on sphere $S_{l}$ for $l \neq j$.

\begin{lemma}
Let the outgoing expansion from sphere $S_j$ be given by
\begin{eqnarray*} 
\bE^j_0 &=&  \sum_{n,m} a_{n,m}^{j} 
\nabla \times \nabla \times (\xb_j \phi_{n,m}^{k_0}(\xb_j) ) 
\, + \, i \omega \mu_0  \sum_{n,m} b_{n,m}^{j} \nabla \times (\xb_j \phi_{n,m}^{k_0}(\xb_j) ) \\
\bH_0^j &=& \sum_{n,m} b_{n,m}^{j} 
\nabla \times \nabla \times (\xb_j \phi_{n,m}^{k_0}(\xb_j) ) 
\, - \, i \omega \epsilon_0 \sum_{n,m} a_{n,m}^{j} \nabla \times (\xb_j \phi_{n,m}^{k_0}(\xb_j) ) .
\end{eqnarray*}
Then, the corresponding field induced on the surface of sphere $S_l$ can be represented in the form
\begin{eqnarray*} 
\bE^l_0 &=&  \sum_{n,m} \gamma_{n,m}^{j,l} 
\nabla \times \nabla \times (\xb_l \psi_{n,m}^{k_0}(\xb_l) ) 
\, + \, i \omega \mu_0  \sum_{n,m} \delta_{n,m}^{j,l} \nabla \times (\xb_l \psi_{n,m}^{k_0}(\xb_l) ) \\
\bH_0^l &=& \sum_{n,m} \delta_{n,m}^{j,l} 
\nabla \times \nabla \times (\xb_l \psi_{n,m}^{k_0}(\xb_l) ) 
\, - \, i \omega \epsilon_0 \sum_{n,m} \gamma_{n,m}^{j,l} \nabla \times (\xb_l \psi_{n,m}^{k_0}(\xb_l) ) .
\end{eqnarray*}
We denote the mappings from the $\{a_{n,m}^j \}$ and $\{b_{n,m}^j\}$ coefficients to the $\{\gamma_{n,m}^{j,l}\}$ and 
$\{\delta_{n,m}^{j,l}\}$ coefficients by $T^{a,\gamma}_{j,l}$, $T^{b,\gamma}_{j,l}$, 
$T^{a,\delta}_{j,l}$, and $T^{b,\delta}_{j,l}$,
respectively. Each of these mappings depends on the vector from the center of sphere $S_j$ to sphere $S_l$
and the parameters $(\mu_0,\epsilon_0,\omega)$.
\end{lemma}

For convenience, we will sometimes denote vectors of coefficients such as
$\{a_{n,m}^j \}$ by $\vec{a^j}$.
The individual components of a translated vector such as
$T^{a,\delta}_{j,l} \vec{\delta^j}$ will be denoted by $[T^{a,\delta}_{j,l} \vec{\delta^j}]_{n,m}$.

\begin{remark}
The formulae for the {\em translation operators} $T_{j,l}^{a,\gamma}$, $T_{j,l}^{b,\gamma}$, $T_{j,l}^{a,\delta}$, and $T_{j,l}^{b,\delta}$ are rather involved \cite{Epton,Gumerov,MTL}.
If the expansions are truncated at $n=p$ terms, there are $2(p+1)^2$ nonzero coefficients
in both the outgoing $(a_{n,m}^j,b_{n,m}^j)$ and incoming $(\gamma_{n,m}^{j,l},\delta_{n,m}^{j,l})$ 
representations. Each translation operator is dense and, therefore requires $O(p^4)$ operations to apply.
More efficient schemes  \cite{FMMsimple,Gumerov} reduces the cost to $O(p^3)$, while the diagonal-form
of the FMM \cite{CGR06,HFFMM} reduces the cost to $O(p^2 \log p)$ for well-separated spheres in the high-frequency regime. 
\end{remark}

Let us now assume that all outgoing and incoming expansion are truncated at $n=p$ terms. The choice
of $p$ is determined by accuracy considerations. It must be sufficiently large to resolve the $\bE$ and
$\bH$ fields on each sphere surface to the desired precision.
 
Using the preceding lemma, the total field immediately exterior to sphere $S_l$ can be written in the form
\begin{eqnarray} 
\bE^l_0 &=& \bE_l^{in}  \, + \, \sum_{\substack{j=1 \\ j \neq l}}^M 
[ T_{j,l}^{a,\gamma} \vec{a^j} + T_{j,l}^{b,\gamma} \vec{b^j}  ]_{n,m}  \,
\nabla \times \nabla \times (\xb_l \psi_{n,m}^{k_0}(\xb_l) ) \nonumber \\
&& + \, i \omega \mu_0  \sum_{j=1 \atop j \neq l}^M 
[ T_{j,l}^{a,\delta} \vec{a^j} + T_{j,l}^{b,\delta} \vec{b^j}  ]_{n,m} \,
\nabla \times (\xb_l \psi_{n,m}^{k_0}(\xb_l)  ) \label{eq:El0}  \\
 && +  \sum^p_{n,m} a_{n,m}^{l} 
\nabla \times \nabla \times (\xb_l \phi_{n,m}^{k_0}(\xb_l) ) 
\, + \, i \omega \mu_0   \sum^p_{n,m} b_{n,m}^{l} \nabla \times (\xb_l \phi_{n,m}^{k_0}(\xb_l)  ) 
\nonumber \\
\bH^l_0 &=& \bH_l^{in} \, + \, \sum_{j=1 \atop j\neq l}^M   
[ T_{j,l}^{a,\delta} \vec{a^j} + T_{j,l}^{b,\delta} \vec{b^j}  ]_{n,m} \,
\nabla \times \nabla \times (\xb_l \psi_{n,m}^{k_0}(\xb_l)  ) \nonumber \\ 
&& - \, i \omega \epsilon_0 \sum_{j=1 \atop j\neq l}^M 
[ T_{j,l}^{a,\gamma} \vec{a^j} + T_{j,l}^{b,\gamma} \vec{b^j}  ]_{n,m}  \,
\nabla \times (\xb_l \psi_{n,m}^{k_0}(\xb_l)  ) \label{eq:Hl0}  \\
 && +  \sum^p_{n,m} b_{n,m}^{l} 
\nabla \times \nabla \times (\xb_l \phi_{n,m}^{k_0}(\xb_l) ) 
\, - \, i \omega \epsilon_0   \sum^p_{n,m} a_{n,m}^{l} \nabla \times (\xb_l \phi_{n,m}^{k_0}(\xb_l)  ) \, . 
\nonumber
\end{eqnarray}
The first terms in the preceding expressions for $\bE_0^l,\bH_0^l$ account for the incoming field,
while the next two terms account for the scattered field coming from all other spheres.
The last two terms in each expression account for the fields being scattered by $S_l$ itself.

It is now clear how to 
apply the interface conditions (\ref{mwaveint1}). We simply equate
the tangential components of $\bE_0^l,\bH_0^l$ defined in 
(\ref{eq:El0}),(\ref{eq:Hl0}) with the tangential components of 
the interior representations $(\bE_l,\bH_l)$ defined in (\ref{Elint}),(\ref{Hlint}).
This yields a dense linear system of dimension
$4M(p+1)^2$ for the coefficients $(a_{n,m}^l, b_{n,m}^l, c_{n,m}^l, d_{n,m}^l)$.
We will refer to this system as the {\em multiple scattering equations}. Writing the equations out
explicitly is not especially informative, and we omit it. 

\begin{remark}
The scattering matrix ${\cal S}$ (Definition \ref{scatmatdef}) 
allows for the elimination of the interior variables $(c_{n,m}^l, d_{n,m}^l)$, so that the 
one can solve a modified system of dimension $2M(p+1)^2$ for the 
coefficients $(a_{n,m}^l, b_{n,m}^l)$ describing the exterior field alone. 

\begin{eqnarray}
\left( \begin{array}{c}
a^l_{n,m} \\ 
\\
b^l_{n,m} 
\end{array} \right) 
&=& {\cal S}
\left(  \begin{array}{c}  
 \alpha^l_{n,m} + \displaystyle{\sum_{\substack{ j=1 \\ j \neq l}}^M 
 [ T_{j,l}^{a,\gamma} \vec{a^j} + T_{j,l}^{b,\gamma} \vec{b^j} ]_{n,m}} \\ 
 \\
 \beta^l_{n,m} + \displaystyle{\sum_{\substack{j=1 \\ j\neq l}}^M   
 [ T_{j,l}^{a,\delta} \vec{a^j} + T_{j,l}^{b,\delta} \vec{b^j} ]_{n,m}}  \end{array}  \right)
\end{eqnarray}
\end{remark}

It is worth emphasizing that the multiple scattering equations are hardly new. There is a vast
literature on the subject, which we do not seek to review here. We refer the reader to the textbooks
 \cite{Bohren,Born,Hafner,Martin,MTL} and the papers \cite{Gumerov,Xu}.

\subsection{Iterative solution of the multiple scattering problem for a system of spheres}

We will solve the multiple scattering equations iteratively, using GMRES \cite{GMRES} with a 
block diagonal preconditioner, each block corresponding to the unknowns on a single sphere.
In applying the preconditioner, we simply invert each of the $M$ diagonal blocks, 
which corresponds to solving the single sphere scattering problem described in section \ref{sec:single}.
Since all $M$ spheres interact, however, the system matrix is dense.
Each matrix-vector multiply in the iterative solution process, 
if carried out naively, would require $O(M^2 \, p^3)$ work. 

In order to accelerate the solution procedure, the wideband fast multipole method (FMM) \cite{CGR06} 
can easily be modified to reduce the cost to $O(M p^3)$ work per iteration. This is discussed in
the context of  acoustic scattering in \cite{Gumerov2,Hesford}. 
Since the literature on FMMs is substantial, we omit a detailed discussion of the technique, 
but present results in section \ref{sec:results}.

\section{Scattering from an arbitrary inclusion} \label{sec:muller}

\begin{figure}[h]
\centering
\includegraphics[width=1.5in]{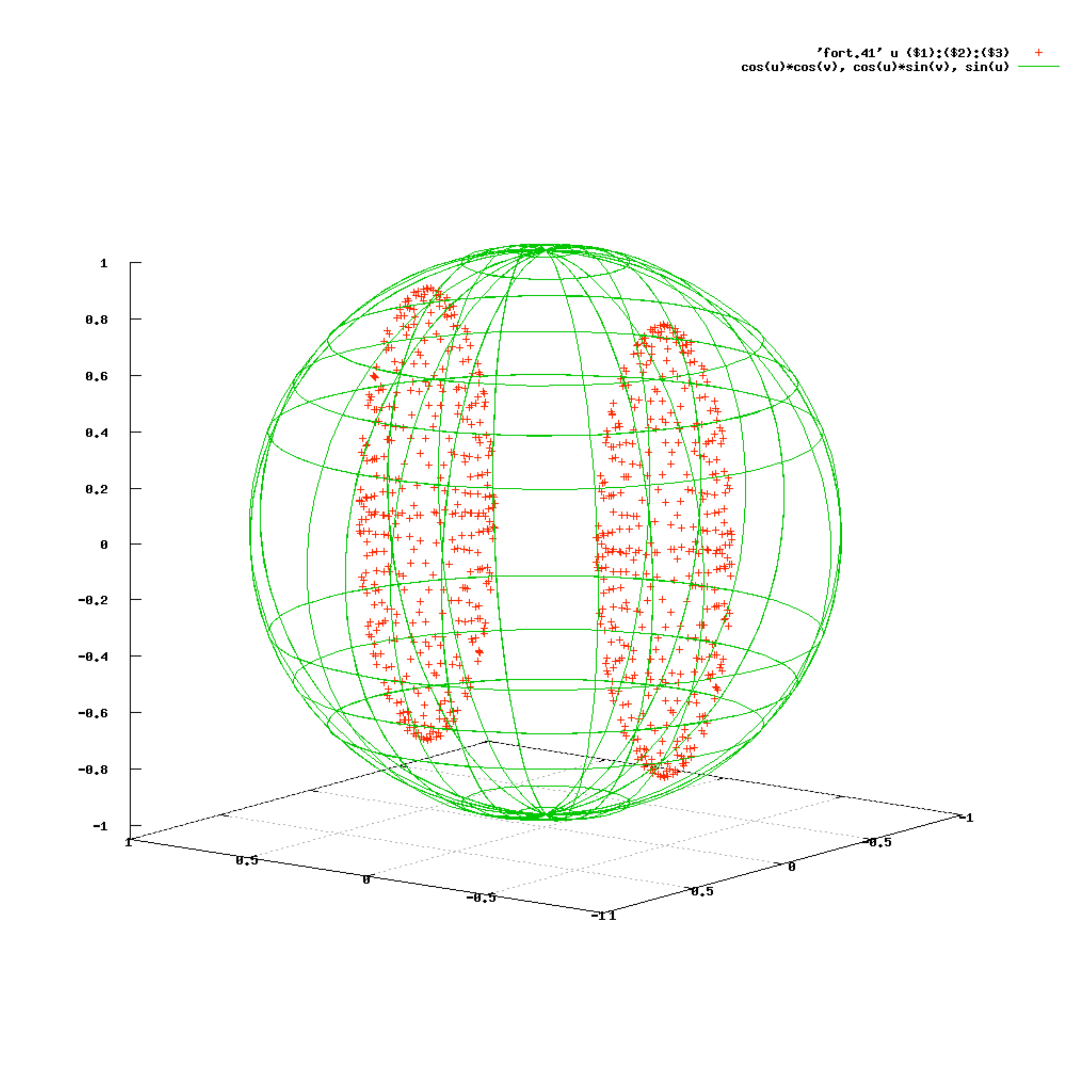}
\caption{\sf A pair of triangulated ellipsoids define a bounded domain $D_1$ that lies with an enclosing sphere $S_1$. The scattering matrix for $D_1$ will be created on $S_1$ and used to represent the 
exterior field.
\label{fig:scatgeom}}
\end{figure}

Suppose now that instead of a sphere, we are given a smooth inclusion (or set of inclusions) 
$D_1$ with 
permittivity $\epsilon_1$ and permeability $\mu_1$ embedded in the same infinite medium as above.
We will suppose further that $D_1$ can be enclosed in a sphere $S_1$
(Fig. \ref{fig:scatgeom}). As before,
at the material interface, the conditions to be satisfied are
(\ref{mwaveint1}). The Debye-Lorenz-Mie formalism cannot be applied in this case, and
attempts to do so (called the T-matrix method) suffer from ill-conditioning when $D_1$ is 
sufficiently non-spherical. We, therefore, turn to the standard representation of electromagnetic
fields in general geometries, based on the vector and scalar potentials and anti-potentials
\cite{Muller,Papas}.

The vector potential in domain $l$ ($l = 0,1$) is defined by 
\[ \A_l(\xb) = \mu_l \int_{\partial D_1}  g_l(\xb-\yb) \, \J_l(\yb) \, ds_\yb \]
and 
$g_l(\xb) = e^{i k_l \|\xb\|}/\|\xb\|$ with 
$k_l = \sqrt{ w^2 \epsilon_j \mu_j}$. When the argument of the square root is 
complex, $k_l$ is taken to lie in the upper half-plane.
We define the vector anti-potential in domain $l$ by 
\[ {\tilde \A}_l(\xb) = \epsilon_l \int_{\partial D_1} 
 g_l(\xb-\yb) \, \K_l \, ds. \]
From these, we may write
\begin{eqnarray*}
\E_l &=& -\nabla \phi_l + i \omega \A_l  - 
       \frac{1}{\epsilon_l} \nabla \times {\tilde \A}_l \\
\H_l &=& \frac{1}{\mu_l} \nabla \times \A_l 
- \nabla \psi_l + i \omega {\tilde \A}_l. 
\end{eqnarray*}
where
\begin{eqnarray*}
\phi_l &=& \frac{1}{i\omega \epsilon_l \mu_l} \nabla \cdot \A_l \\
\psi_l &=& \frac{1}{i\omega \epsilon_l \mu_l} \nabla \cdot {\tilde \A}_l.
\end{eqnarray*}

\begin{figure}[h]
\centering
\vspace{.3in}
\includegraphics[width=3.5in]{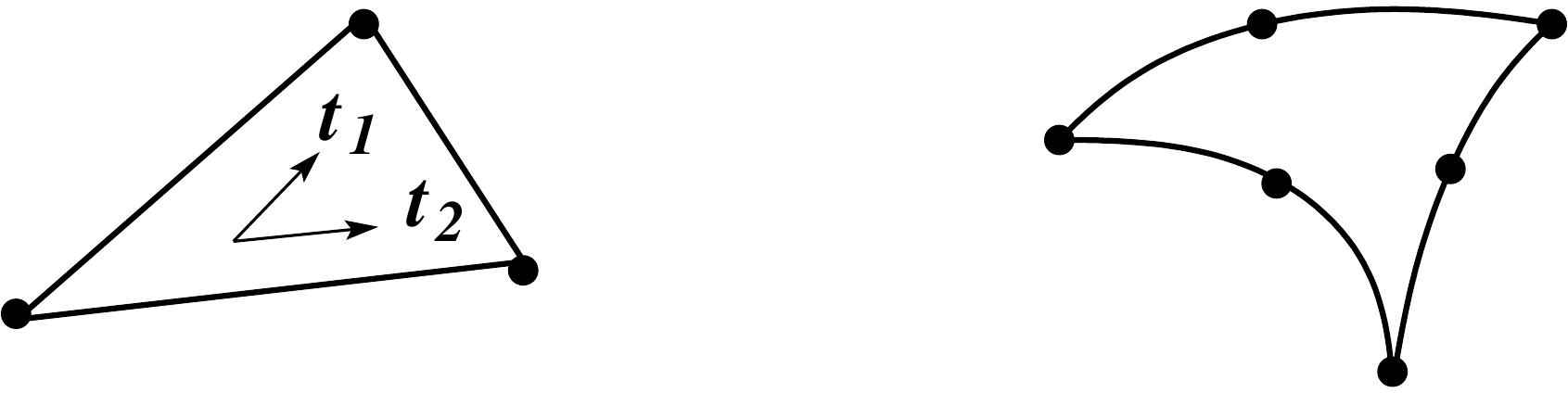}
\caption{\sf In the simplest geometric model, the surface of the scatterer $\partial D_1$ is
approximated by a collection of flat triangles, defined by the locations of its three vertices in $\rthree$.
On each triangle, there are two two linearly independent tangent directions ${\bf t}_1$ and ${\bf t}_2$.
The unknown electric and magnetic currents $\J$ and $\K$ on each triangle are defined by
$j_1 {\bf t}_1 +  j_2 {\bf t}_2$ and $k_1 {\bf t}_1 +  k_2 {\bf t}_2$, respectively, and the electromagnetic
fields are evaluated at the triangle centroids.
For higher order accuracy, each quadratic surface patch is specified by six nodes:
the three triangle vertices and three additional points, one on 
each curved triangle side. Three ``support nodes"  ${\bf x}^1, {\bf x}^2, {\bf x}^3$ 
are then selected in the interior of each patch. Our representation for 
$\J$ and $\K$ at each support node ${\bf x}^i$ is of the form  
$j^i_1 {\bf t}^i_1 +  j^i_2 {\bf t}^i_2$ and $k^i_1 {\bf t}^i_1 +  k^i_2 {\bf t}^i_2$,
where ${\bf t}^i_1,{\bf t}^i_2$ are linearly independent tangent vectors at ${\bf x}^i$. The 
support nodes are also the points where we evaluate the electromagnetic fields and impose
interface conditions.
\label{fig:curvetri}}
\end{figure}

As written above, we have twelve degrees of freedom at each point 
$P \in \partial D_1$,
namely the three Cartesian components of $\J_0,\J_1, \K_0,
\K_1$, but only four boundary conditions (the continuity of the tangential components of 
$\bE$ and $\bH$).
We will assume, however, that the functions 
$\J_0,\J_1, \K_0, \K_1$ are surface currents and that
the following linear relations hold
\[  \J_0 =   \frac{\epsilon_0}{\epsilon_1}  \J_1 \qquad 
\K_0 = \frac{\mu_0}{\mu_1} \K_1  \, .
 \]
This leaves four degrees of freedom. Imposing the conditions  (\ref{mwaveint1})
on $\J_1,\K_1$ results in M\"{u}ller's integral equation \cite{Muller}, a resonance-free Fredholm
equation of the second kind.

In more detail, using the facts that
\[ \nabla_{\xb} \times (g_l(\xb - \yb) \, \K(\yb)) = \nabla_{\xb} 
g_l \times \K(\yb) \, , \]
\[ a \times b \times c = b(a \cdot c) - c(a\cdot b) \, , \]
and, for $\by_0 \in \partial D_1$,
\[ \lim_{\xb \rightarrow \yb_0 \atop \xb \in D_0}
\int_{\partial D_1} \frac{\partial g_l}{\partial n_{\yb_0}}(\xb-\yb) \, \sigma(\yb) ds_\yb = 
\frac{1}{2} \sigma(\yb_0) + 
\oint_{\partial D_1}  \frac{\partial g_l}{\partial n_{\yb_0}}(\yb_0-\yb) \, \sigma(\yb) ds_\yb \]
\[ \lim_{\xb \rightarrow \yb_0 \atop \xb \in D_1}
\int_{\partial D_1} \frac{\partial g_l}{\partial n_{\yb_0}}(\xb-\yb) \, \sigma(\yb) ds_\yb = 
-\frac{1}{2} \sigma(\yb_0) + 
\oint_{\partial D_1}  \frac{\partial g_l}{\partial n_{\yb_0}}(\yb_0-\yb) \, \sigma(\yb) ds_\yb, \]
we obtain the following coupled set of equations.

\begin{eqnarray}
-\n \times \E^{in} &=& \frac{i \omega}{\epsilon_1} 
\int_{\partial D_1} [ \epsilon_0 \mu_0  \, g_0 - 
\epsilon_1\mu_1 \, g_1] \, (\n \times \J_1) \, ds_\yb \nonumber \\
& & + 
\frac{i}{\omega \epsilon_1} \n \times \int_{\partial D_1} 
[ \nabla \nabla g_0  -  \nabla \nabla g_1]  \, \J_1  \, ds_\yb
\label{cancelE} \\
& & - \mu_0 \int_{\partial D_1} 
\left( \frac{\nabla g_0}{\mu_1 \epsilon_0} - 
\frac{\nabla g_1}{\mu_0 \epsilon_1} \right) \, (\n \cdot \K_1)  \, ds_\yb 
\nonumber \\
& & + \left( \frac{1}{2\epsilon_0} + \frac{1}{2\epsilon_1} \right) \K_1
+ \mu_0  \oint_{\partial D_1}
\left(
\frac{1}{\mu_1 \epsilon_0}
\frac{\partial g_0}{\partial n} -
\frac{1}{\mu_0 \epsilon_1}
\frac{\partial g_1}{\partial n} \right) \, \K_1 \, ds_\yb
\nonumber
\end{eqnarray}

\begin{eqnarray}
- \n \times \H^{in} &=& \frac{i \omega}{\mu_1}
\int_{\partial D_1} [ \mu_0 \epsilon_0 \, g_0 - \epsilon_1\mu_1 \, g_1] 
\, (\n \times \K_1)\, ds_\yb \nonumber \\
& & + 
\frac{i}{\omega\mu_1}  \n \times \int_{\partial D_1} 
[ \nabla \nabla g_0   - \nabla \nabla g_1] \,  \K_1 \, ds_\yb 
\label{cancelH} \\
& & - \epsilon_0 \int_{\partial D_1} 
\left( \frac{ \nabla g_0}{\epsilon_1 \mu_0} - 
\frac{\nabla g_1}{\epsilon_0 \mu_1} \right) \, (\n \cdot \J_1)  \, ds_\yb \nonumber \\
& & + \left( \frac{1}{2\mu_0} + \frac{1}{2\mu_1} \right) \J_1
+ \epsilon_0 \oint_{\partial D_1} 
\left( \frac{1}{\epsilon_1 \mu_0}
\frac{\partial g_0}{\partial n} -
\frac{1}{\epsilon_0 \mu_1}
\frac{\partial g_1}{\partial n} \right)
\, \J_1 \, ds_\yb .
\nonumber
\end{eqnarray}

Because the M\"{ u}ller equation is a second kind Fredholm equation, the order of accuracy 
of the solution is that of the underlying quadrature rule. For first order accuracy, we assume
$\J_1$ and $\K_1$ are piecewise constant current densities on a flat triangulated surface.
For second order accuracy, we assume $\J_1$ and $\K_1$ are piecewise linear current densities 
on a piecewise quadratic surface with each curved triangle defined by six points (Fig. \ref{fig:curvetri}).

For each discretization node, we evaluate the relevant electromagnetic field component using
a mixture of analytic and numerical quadratures on each triangle. More precisely, we use the
method of singularity subtraction - computing integrals analytically for the kernel $1/r$ and its derivatives
and using numerical quadrature for the {\em difference kernel} $[e^{ikr}-1]/r$, which is smoother.
This results
in a complex linear system of dimension $4N \times 4N$ for first order accuracy and of dimension
$12N \times 12N$ for second order accuracy, where $N$ denotes the number of triangles.
For small $N$, say $N < 1000$, one can use direct LU-factorization to solve the linear system.
For larger values of $N$, iterative solution with FMM-acceleration becomes much more practical
\cite{CGR06,FMMsimple}.

\subsection{The scattering matrix for $D_1$}

Suppose now that we are interested in scattering from the two ellipsoids $D_1$ shown in 
Fig. \ref{fig:scatgeom} due to an incoming field which is regular in the enclosing sphere 
$S_1$. Such an incoming field can be expanded within $S_1$ in the form
\begin{eqnarray*} 
\bE^{in}(\xb) &=&  \sum_{n,m} \alpha_{n,m}
\nabla \times \nabla \times (\xb \psi_{n,m}^{k_0}) 
\, + \, i \omega \mu_0 \sum_{n,m} \beta_{n,m} \nabla \times (\xb \psi_{n,m}^{k_0}) \\
\bH^{in}(\xb) &=&  \sum_{n,m} \beta_{n,m}
\nabla \times \nabla \times (\xb \psi_{n,m}^{k_0}) 
\, - \, i \omega \epsilon_0 \sum_{n,m} \alpha_{n,m} \nabla \times (\xb \psi_{n,m}^{k_0}) \, ,
\end{eqnarray*}
as in Section \ref{sec:single}. Each  (vector) spherical harmonic
modes, corresponding to a single $\alpha_{n,m}$ or $\beta_{n,m}$, defines a particular
incoming field on $D_1$. More precisely, we can solve the M\"{ u}ller equation for a right-hand 
side obtained by setting the incoming field to be
\begin{equation}
 \bE^{in}_{1,n,m}(\xb) = \, \nabla \times \nabla \times (\xb \psi_{n,m}^{k_0}),  \quad
\bH^{in}_{1,n,m}(\xb) = \, - \, i \omega \epsilon_0 \nabla \times (\xb \psi_{n,m}^{k_0}) \, ,
\end{equation}
corresponding to setting a fixed $\alpha_{n,m} = 1$ and all other coefficients to zero. 
Similarly, we can set the incoming field to be  
\begin{equation}
\bE^{in}_{2,n,m}(\xb) =  \, + \, i \omega \mu_0  \nabla \times (\xb \psi_{n,m}^{k_0}), \quad
\bH^{in}_{2,n,m}(\xb) = \,  \nabla \times \nabla \times (\xb \psi_{n,m}^{k_0})
\end{equation}
corresponding to setting a fixed $\beta_{n,m} = 1$ and all other coefficients to zero. 
We can then store either the electric and magnetic currents
$\J^{1,n,m}_1, \K^{1,n,m}_1$ or 
$\J^{2,n,m}_1, \K^{2,n,m}_1$ 
induced by these (unit) incoming fields or just convert these currents to the coefficients of the outgoing
(scattered) fields:
\[
\bE^{sc}_{1,n,m}(\xb) =  \sum_{n',m'} a^{1,n,m}_{n',m'} 
\nabla \times \nabla \times (\xb \phi_{n',m'}^{k_0}) 
\, + \, i \omega \mu_0 \sum_{n',m'} b^{1,n,m}_{n',m'} \nabla \times (\xb \phi_{n',m'}^{k_0})
\]
\[
\bH^{sc}_{1,n,m}(\xb) =  \sum_{n',m'} b^{1,n,m}_{n',m'} 
\nabla \times \nabla \times (\xb \phi_{n',m'}^{k_0}) 
\, - \, i \omega \epsilon_0 \sum_{n',m'} a^{1,n,m}_{n',m'} \nabla \times (\xb \phi_{n',m'}^{k_0}) 
\] 
and
\[
\bE^{sc}_{2,n,m}(\xb) =  \sum_{n',m'} a^{2,n,m}_{n',m'} 
\nabla \times \nabla \times (\xb \phi_{n',m'}^{k_0}) 
\, + \, i \omega \mu_0 \sum_{n',m'} b^{2,n,m}_{n',m'} \nabla \times (\xb \phi_{n',m'}^{k_0}) \]
\[
\bH^{sc}_{2,n,m}(\xb) =  \sum_{n',m'} b^{2,n,m}_{n',m'} 
\nabla \times \nabla \times (\xb \phi_{n',m'}^{k_0}) 
\, - \, i \omega \epsilon_0 \sum_{n',m'} a^{2,n,m}_{n',m'} \nabla \times (\xb \phi_{n',m'}^{k_0}). 
\]
The formula for converting the currents $\J^{1,n,m}_1, \K^{1,n,m}_1$ 
to the coefficients can be obtained by orthogonal projection of the induced field on the enclosing sphere  
\cite{Papas}.

By superposition, an incoming field defined by the vector of incoming coefficients
$\{\alpha_{n,m}, \beta_{n,m} \}$ results in a scattered field of the form
\begin{eqnarray*} 
\bE^{sc}(\xb) &=&  \sum_{n',m'} a_{n',m'}
\nabla \times \nabla \times (\xb \phi_{n',m'}^{k_0}) 
\, + \, i \omega \mu_0 \sum_{n',m'} b_{n',m'} \nabla \times (\xb \phi_{n',m'}^{k_0}) \\
\bH^{sc}(\xb) &=&  \sum_{n',m'} b_{n',m'}
\nabla \times \nabla \times (\xb \phi_{n',m'}^{k_0}) 
\, - \, i \omega \epsilon_0 \sum_{n',m'} a_{n',m'} \nabla \times (\xb \phi_{n',m'}^{k_0}) \, ,
\end{eqnarray*}
with the coefficients of the scattered field given by 
\begin{eqnarray*}
 a_{n',m'} &=&   \sum_{n,m} \alpha_{n,m} a^{1,n,m}_{n',m'} + \beta_{n,m} a^{2,n,m}_{n',m'} \\
 b_{n',m'} &=&   \sum_{n,m} \alpha_{n,m} b^{1,n,m}_{n',m'} + \beta_{n,m} b^{2,n,m}_{n',m'}.
\end{eqnarray*}
The matrix mapping the incoming to the scattered coefficients is referred to as
the {\em scattering matrix} for the structure $D_1$. 

Fixing the order of the expansions above at $p$,
there are $4p^2$ possible basis functions that span the space of all possible incoming fields.
We must, therefore solve $4 p^2$ M\"{u}ller integral equations on the detailed geometry
defining $D_1$. To store the currents induced by each incoming field requires $O(N p^2)$ memory, 
where $N$ denotes the number of degrees of freedom used in the discretization of the integral equation. 
The scattering matrix itself requires storing $O(16 p^4)$ complex numbers. 
While somewhat expensive, this is a pre-computation step, in anticipation of simulating microstructures
with thousands or millions of inclusions of the same identical shape, but well enough separated
that the scattering matrices are accurate. 

\section{Multiple scattering from well-separated non-spherical inclusions}

Once the scattering matrix is known, the solution to the full Maxwell equations for geometries with $N$
inclusions ($N = 200$ in Fig. \ref{f:1})  can be turned into a multiple-scattering
problem based only on the enclosing spheres. 
That is, the inclusions can be replaced by their scattering matrices and the multiple scattering 
method of section \ref{multspherescat} can be used with trivial modifications. 

There are two distinct advantages to be gained here. First,
we have reduced the number of degrees of freedom from, say, 5,000 or 10,000 unknown 
current density values per inclusion to, say, 400 expansion coefficients. Just as 
important, however, is that we have precomputed
the solution operator for each inclusion in isolation, so that the linear system we solve
by iteration on the multi-sphere system
is much more well-conditioned. Further, the FMM reduces the cost of each iteration 
from $O(N^2)$ to $O(N \log N)$ and is particularly efficient here, since the complicated quadratures
on triangulated surfaces have been subsumed into the precomputation step.

The principal limitations of the method are 1) that  some modest separation
distance between inclusions is required and 2)
that the bookkeeping  becomes a bit awkward if more than a few distinct nanoparticle types are allowed. 
In many experimental settings, both conditions are satisfied. 

It is worth noting that the method of this paper can be viewed as a 
{\em reduced order model} for the scattering problem. In broad terms, 
the idea is not new and there is substantial activity in this area in both electromagnetics 
and other fields (see, for example, \cite{Hesthaven}).
It is also worth noting that the method is ``rigorous" in the sense that the error is determined in a
straightforward manner by the accuracy of the M\"{ u}ller integral equation solver and the order of expansion
of the scattering matrix. It fails (or needs local modification) if and only if two enclosing spheres intersect.

\section{Numerical Examples} \label{sec:results}

As discussed in section \ref{sec:muller}, the M\"{ u}ller integral equation is an effective
method for determining the scattering matrix from a dielectric inclusion of arbitrary shape. 
To illustrate its performance, we consider the geometry in Fig. \ref{fig:scatgeom}, consisting 
of a pair of ellipsoids triangulated with piecewise quadratic triangles on which we allow piecewise
linear current densities. 
Each triangle has three nodes with two degrees of freedom for each current
(electric and magnetic) at each point, resulting in a complex linear system of dimension
$2160 \times 2160$. 
(All calculations and timings reported in this section  have been carried out
using a 12-core 2.93GHz Intel Xeon workstation.)
LU factorization requires 3.5 seconds, and the subsequent solution
requires 0.1 seconds for each possible incoming mode. 
With 720 triangles, the linear system has dimension 8640 and with
2880 triangles, there are 34,560 degrees of freedom. These require 106 and 2,620 seconds 
to factor, respectively. The solution times for each incoming mode are 3.52 and 87 seconds, 
respectively. We could accelerate these solution times using fast multipole-based codes 
(or any of a variety of other ``fast" algorithms), but we view
this cost as an initialization step and the CPU times are acceptable. The errors are of the order
$10^{-3}$, $10^{-4}$, and $10^{-5}$ for the successively finer discretizations, somewhat better
than the expected second order convergence.

To illustrate the performance of the FMPS algorithm, we consider a $21 \times 21 \times 2$ array of
scatterers, each consisting of an ellipsoid pair with a scattering matrix derived from the 
M\"{ u}ller integral equation of order $p=3$. Using the same 12-core 2.93GHz Intel Xeon workstation,
the time required was about 2 seconds per iteration, with six iterations required for GMRES to converge to 3
digits.  The ``slow" multiple scattering (SMPS) approach, without fast multipole acceleration, required
about 7 seconds per iteration.
For a $21 \times 21 \times 4$ array, the cost was about 6 seconds per iteration (28 seconds for SMPS) and for 
a $21 \times 21 \times 8$ array, the cost was about 23 seconds per iteration (108 seconds for SMPS).
For a $21 \times 21 \times 16$ array (14,112 ellipsoid pairs), the cost was about 59 seconds per iteration
(440 seconds for SMPS). 

The reason for the modest speedup of the FMPS over the SMPS approach is that the number of spheres
is still rather small. For one million scatterers, the speedup factor would be about 1000. Careful readers
may note that the FMPS scaling appears worse than O(N log N) in successively doubling the simulation from a
 $21 \times 21 \times 4$ array to a $21 \times 21 \times 8$ array to a $21 \times 21 \times 16$ array.
 For those familiar with the FMM, the short explanation is that the  ``near neighbor" cost is not yet in the 
 asymptotic regime in the first two cases. Timings extrapolated from the last case are accurate for any
 volume-filling distribution.

\begin{figure}
\centering
\includegraphics[width=1.52in]{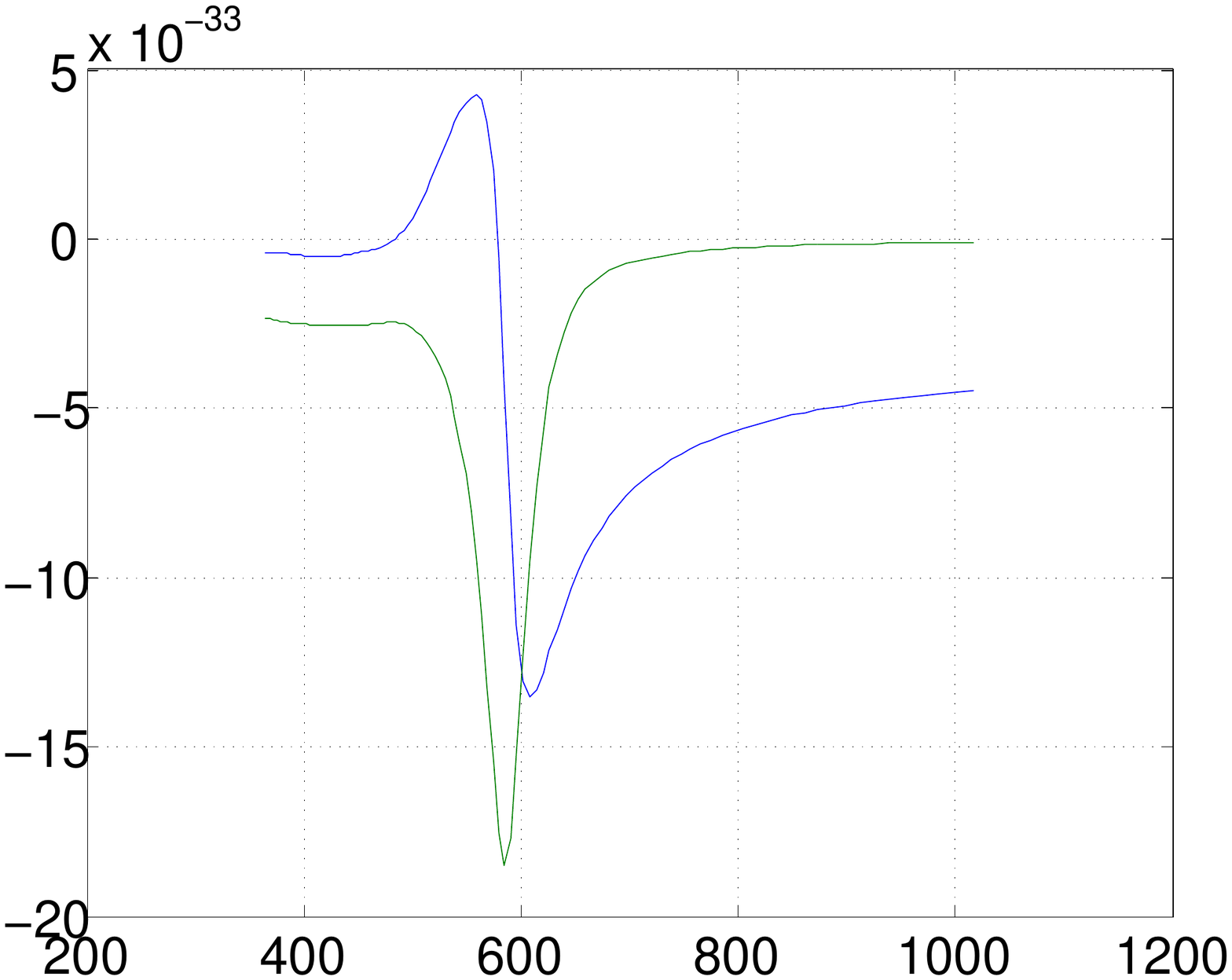} 
\includegraphics[width=1.52in]{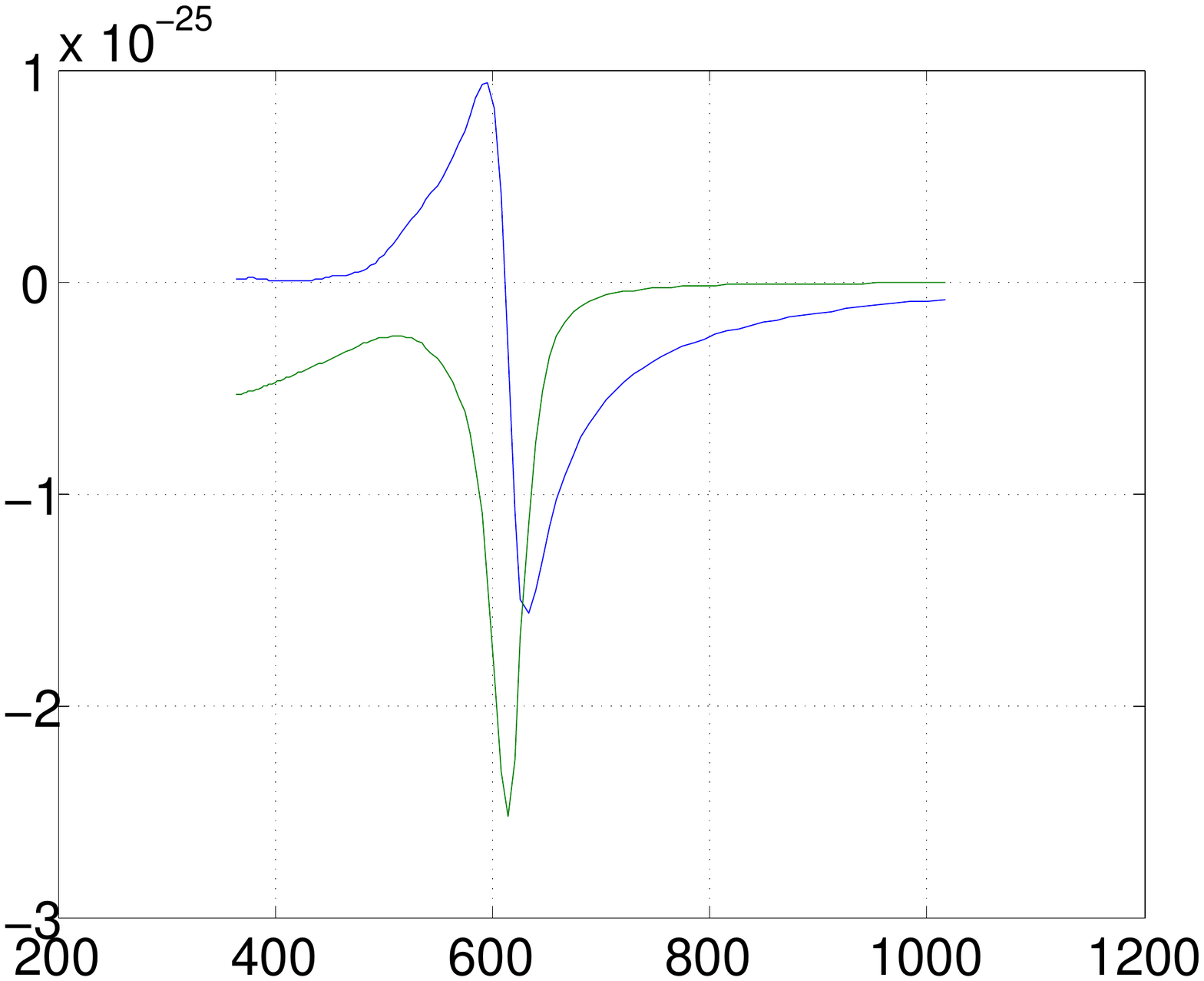} 
\includegraphics[width=1.52in]{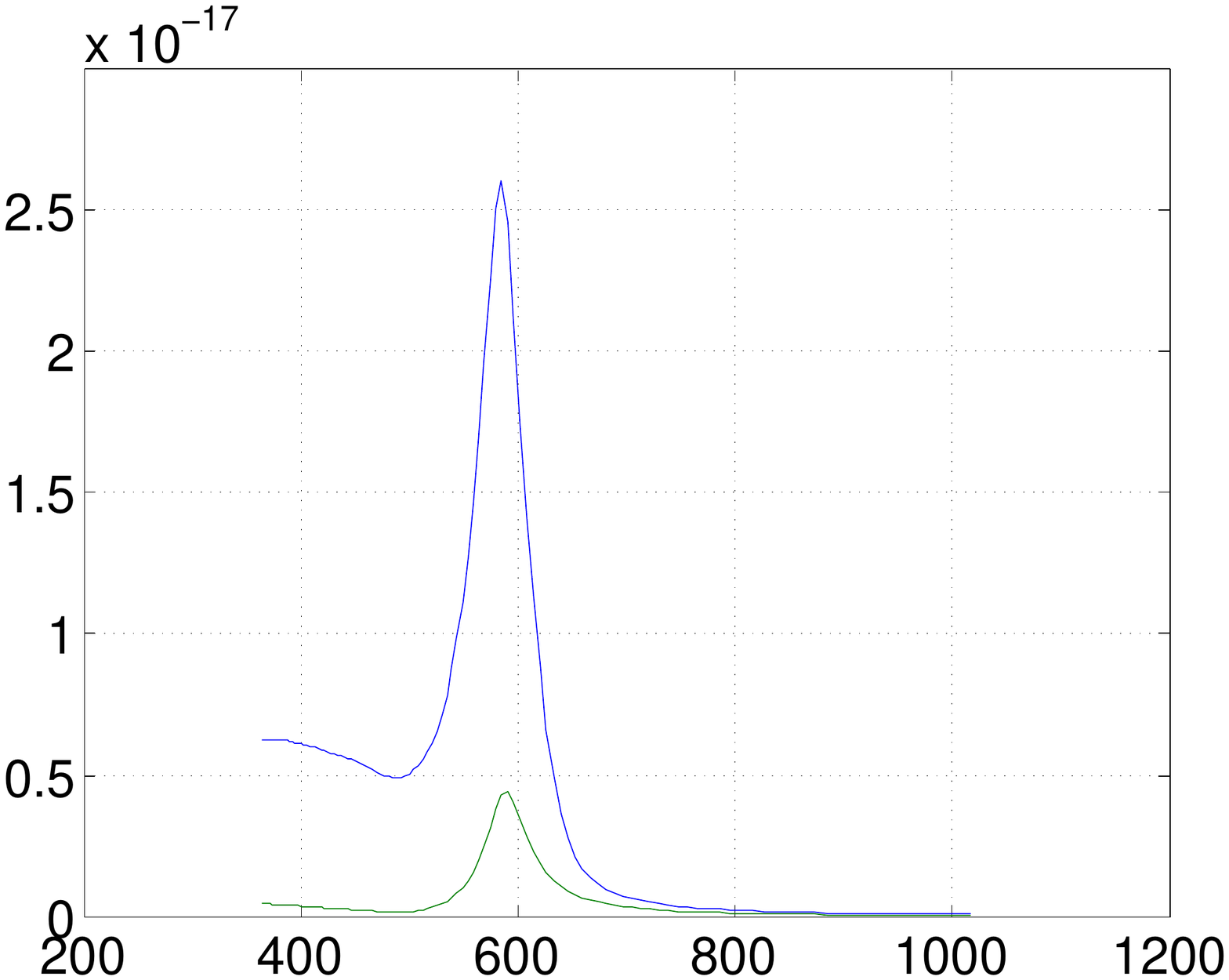} 
\includegraphics[width=1.52in]{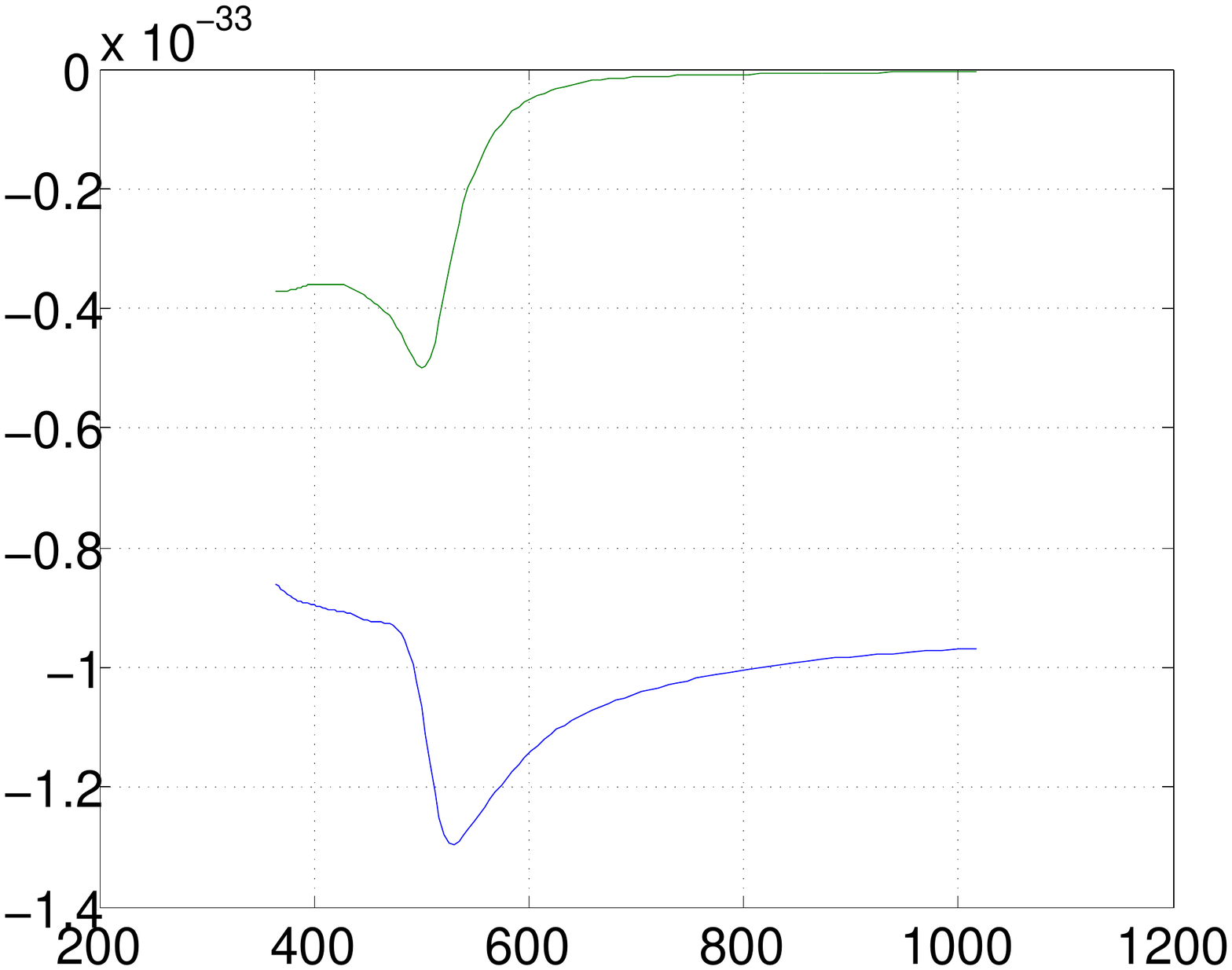} 
\includegraphics[width=1.52in]{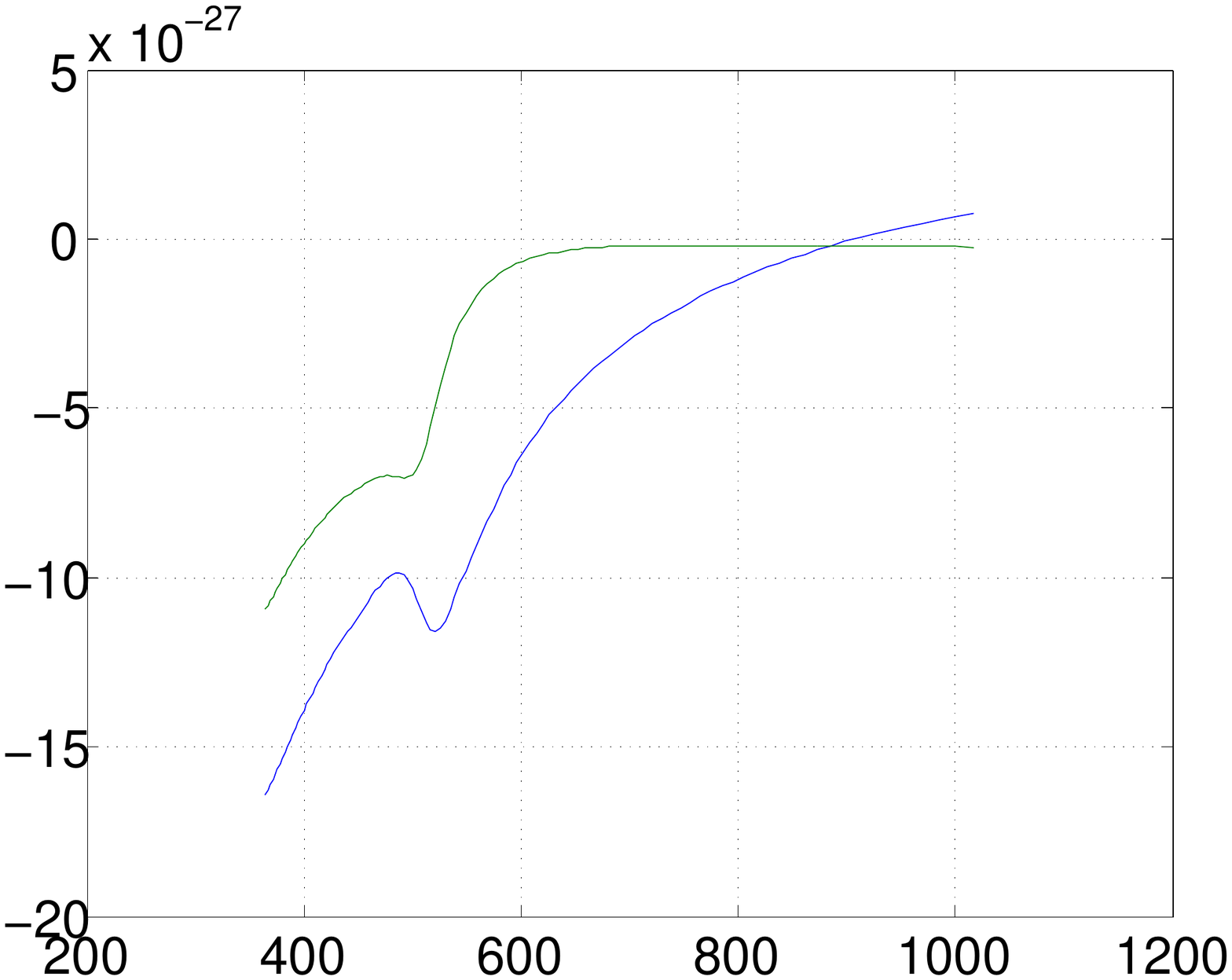} 
\includegraphics[width=1.52in]{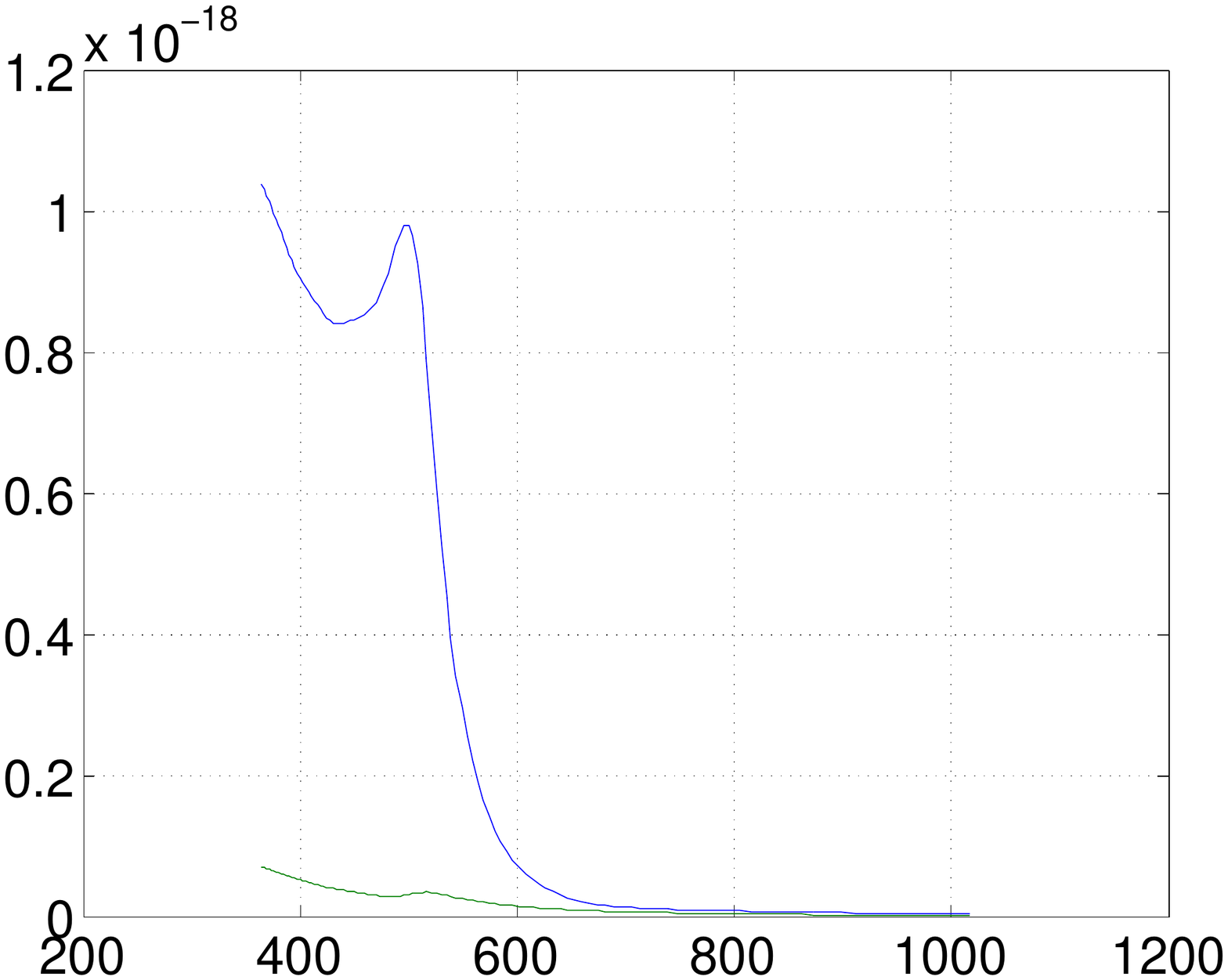} 
\includegraphics[width=1.52in]{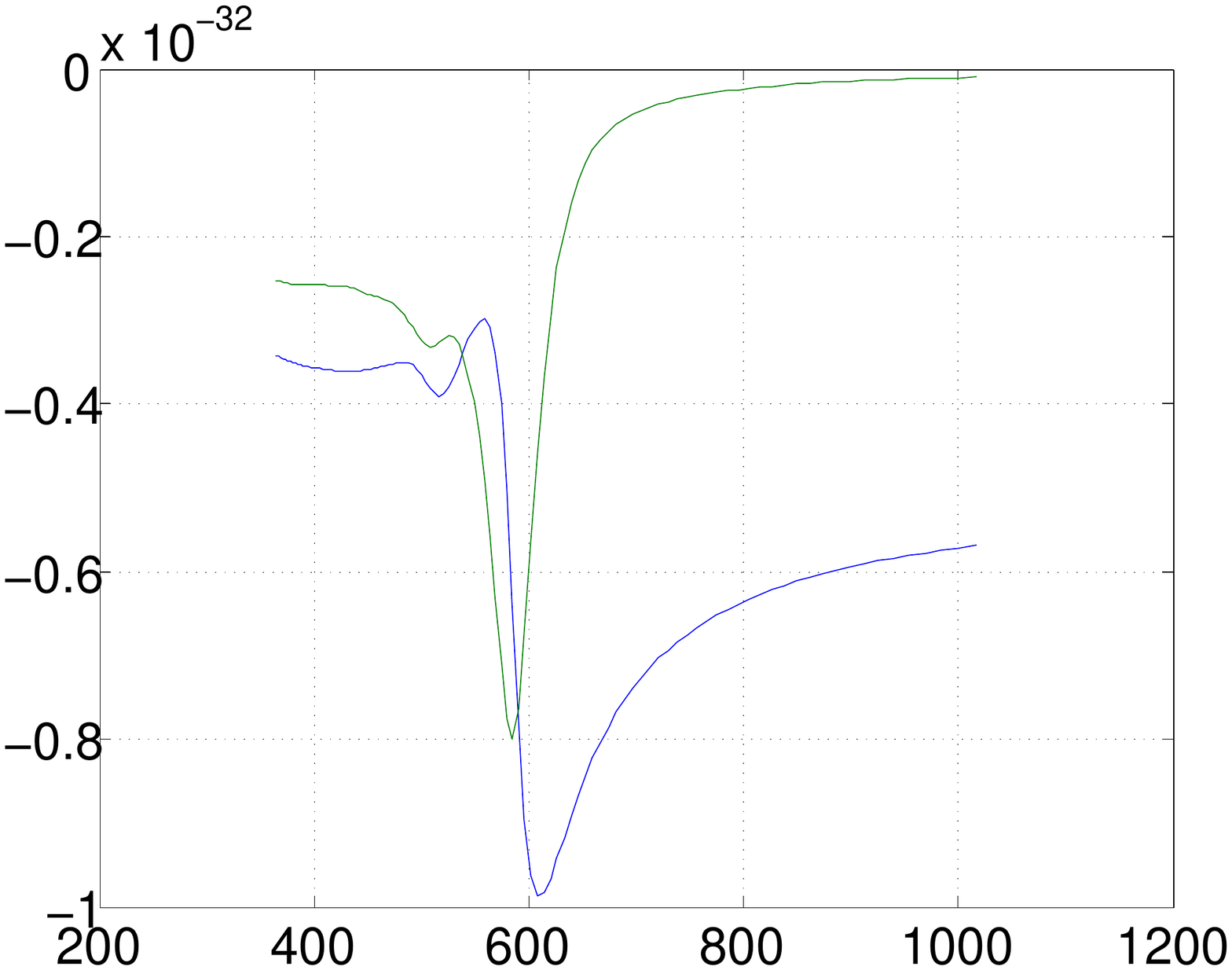} 
\includegraphics[width=1.52in]{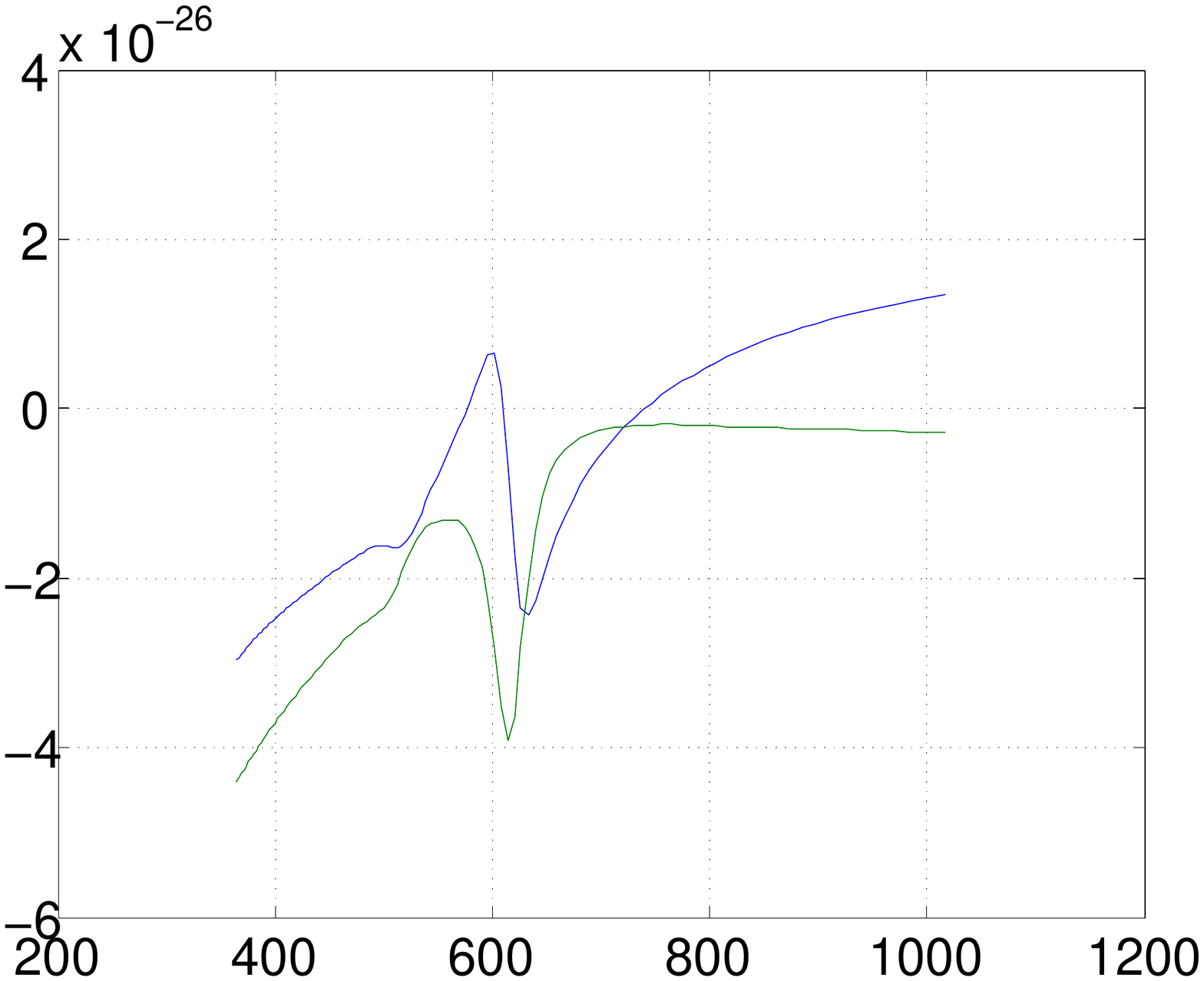} 
\includegraphics[width=1.52in]{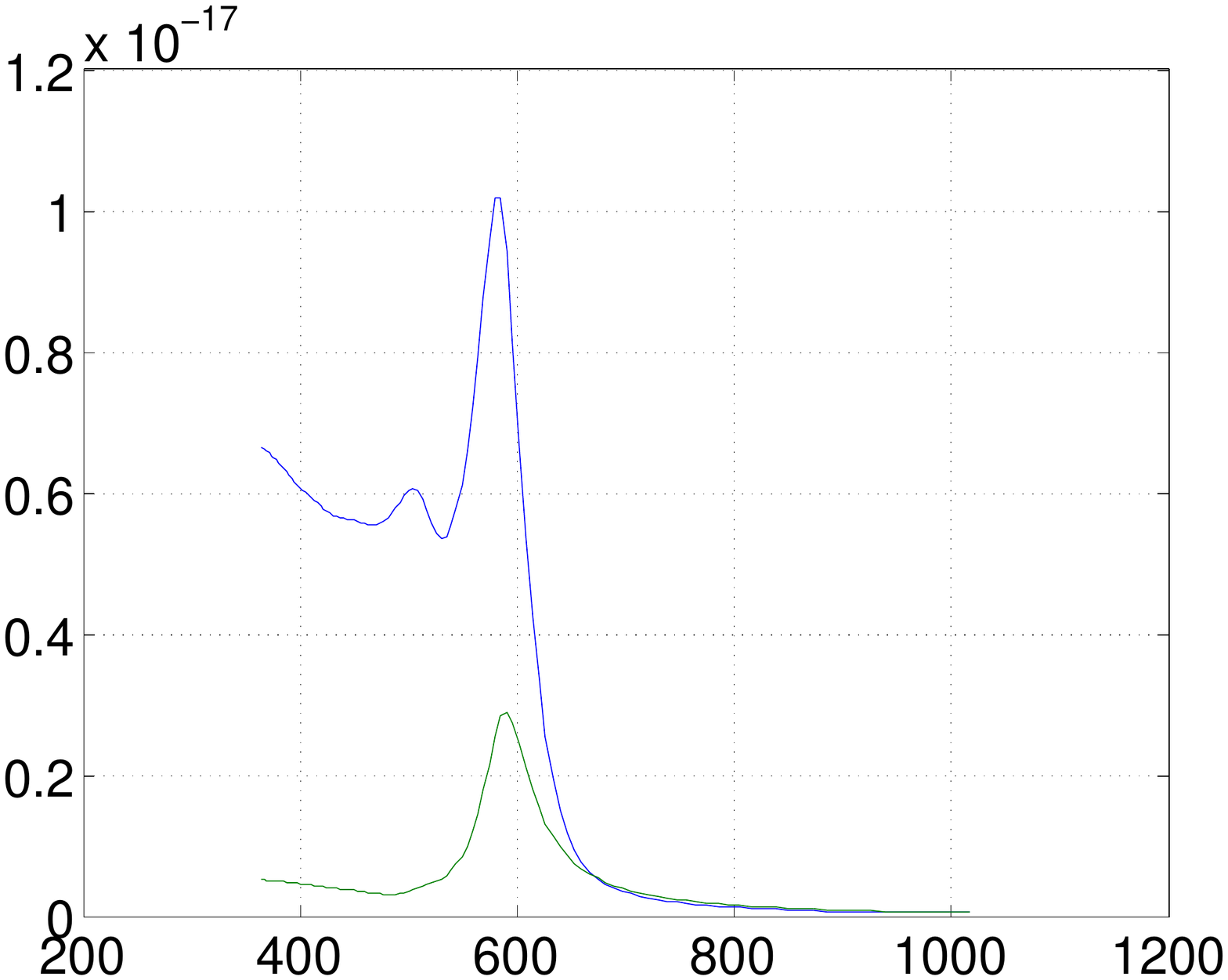} 
\caption{\sf The top row shows a frequency scan of the real and imaginary parts of
electric polarization vector  (left), the real and imaginary parts of
magnetic polarization vector  (middle), and the scattering (right, upper curve) and absorption
(right, lower curve)  for one ellipsoid pair with the long axis oriented parallel to the (linearly polarized)
incoming electric field.
The second row shows the same computed quantities for
one ellipsoid pair with  the long axis oriented parallel to the (linearly polarized)
incoming magnetic field. The third row shows the same computed quantities for four pairs of 
randomly oriented ellipsoid pairs.
\label{f:ex3}}
\end{figure}

Finally, we illustrate the use of the FMPS algorithm in carrying out frequency scans  
for (a) one ellipsoid pair with the long axis oriented parallel to the (linearly polarized)
incoming electric field, (b)one ellipsoid pair with  the long axis oriented parallel to the (linearly polarized)
incoming magnetic field, or (c) four pairs of randomly oriented ellipsoid pairs (Fig. \ref{f:ex3}).

\section{Conclusions} \label{sec:concl}

The method introduced in this paper (fast multi-particle scattering) combines a highly
accurate integral equation solver with multiple scattering theory,
in order to permit the solution to the full Maxwell equations in configurations
typical of engineered composites (metamaterials). We assume that the geometry consists 
of a large number of inclusions embedded in a homogeneous
background. While we have only included a single type of inclusion geometry in our examples above, it is 
clear that the method can easily be applied to permit several such types, so long as there is a 
modest separation between inclusions. 
FMPS is enormously faster than a full FMM-based solver using the full discretization of the geometry.
With 14,112 ellipsoid pairs (the largest example in the preceding section), 
this would require about 30 million degrees of freedom, many minutes
per iteration, and many more iterations.

In its present form, the method cannot be used for tightly packed configurations, which will require more 
elaborate {\em compression} schemes \cite{Fastdirect}. It does, however, permit workstation-based
simulation with millions of inclusions. We are currently working on extending the method so that it
can handle inclusions embedded in a layered medium.

\bibliographystyle{plain}

\end{document}